\documentclass[11pt]{article}
\usepackage{graphics,amssymb,amsmath,fullpage}

\usepackage{todonotes}

\newcommand{\E}{\mathbb{E}}
\newcommand{\R}{\mathbb{R}}
\begin{document}
\title{Sequential escapes: onset of slow domino regime via a saddle connection}

\author{Peter Ashwin\thanks{email: {\tt p.ashwin@exeter.ac.uk}}, Jennifer Creaser and Krasimira Tsaneva-Atanasova\\
Department of Mathematics and 
EPSRC Centre for Predictive Modelling in Healthcare,\\
University of Exeter,\\
Exeter, EX4 4QJ, UK.}

\maketitle

\begin{abstract}
We explore sequential escape behaviour of coupled bistable systems under the influence of stochastic perturbations. We consider transient escapes from a marginally stable ``quiescent'' equilibrium to a more stable ``active'' equilibrium. The presence of coupling introduces dependence between the escape processes: for diffusive coupling there is a strongly coupled limit (fast domino regime) where the escapes are strongly synchronised while for intermediate coupling (slow domino regime) without partially escaped stable states, there is still a delayed effect. These regimes can be associated with bifurcations of equilibria in the low-noise limit. In this paper we consider a localized form of non-diffusive (i.e pulse-like) coupling and find similar changes in the distribution of escape times with coupling strength. However we find transition to a slow domino regime that is not associated with any bifurcations of equilibria. We show that this transition can be understood as a codimension-one saddle connection bifurcation for the low-noise limit. At transition, the most likely escape path from one attractor hits the escape saddle from the basin of another partially escaped attractor. After this bifurcation we find increasing coefficient of variation of the subsequent escape times.
\end{abstract} 

%
%
\section{Introduction}
\label{sec:intro}

Stochastic dynamical systems are often used to describe the behaviour of physical systems where detailed fast chaotic (e.g. thermal) processes are replaced by an idealised stochastic term. The stochastic perturbations (noise) can lead to escape of the system from a state that is stable at a distribution of times that depends on the deterministic dynamics and the amplitude and nature of the noise. This approach was pioneered in the theromodynamics of chemical reactions by Eyring and Kramers \cite{kramers40brown} and has since found applications in a wide range of physical, biological, medical, chemical and engineering applications: see for example references in \cite{AshCreTsa2017,biorev,physrev}. In particular, the powerful theory of escape times of Friedlin and Wentzell \cite{FriedlinWentzell1998} quantifies asymptotical properties of these escapes using large deviations of stochastic processes in the limit of low noise, see also \cite{Berglund2013,BG2006}. 

If we focus on problems where the escape is a transient, rather than a recurrent process (such as for cell-line differentiation in developmental biology \cite{Wangetal2010}) we can think of the escape process as irreversible and ignore returns. This is the case for stochastically perturbed bistable systems in the case where we start in a ``shallow'' (we call quiescent) attractor and escape to a ``deep'' (we call active) attractor where the potential barrier to escape from the quiescent state is much less than for the active state (in fact, the barrier height differences can be quite small and still lead to vast difference in mean escape time for the low noise limit). We will work in this regime where return escape rates are extremely slow and can be ignored.

An extensive literature has considered stochastic resonance \cite{BG2006} between the escape timescales and other timescales within e.g. forcing of the system. Many authors have looked at the problem of coupled escapes, for example \cite{BFG2007a,BFG2007b,Frankowiczetal1982,Malchow_etal_1983,Neiman1994} and in particular in the case of local non-diffusive (i.e. pulse-like) coupling \cite{Zhang1998,Wood2006,Assis2009}. 
Pulse-coupled systems are characterised by exchanging localised in phase space pulse-like signals. Such systems have found wide applications in biology, e.g. populations of synaptically coupled neurones, flashing fireflies, claps of applauding audiences \cite{Mirollo1990,Ermentrout1990,Stankovski2017} as well as engineering, e.g. wireless sensor networks, impulsive control, swarm robotics, smart materials \cite{Zou2015,Proskurnikov2017,Wang2017,Li2018}. Their importance has motivated numerous studies predominantly focussing on understanding the basic synchronisation properties of networks of pulse-coupled oscillators \cite{Mirollo1990,Ermentrout1990,Stankovski2017,Zou2015,Proskurnikov2017,Wang2017,Mirollo2017}.
Here we consider a state-triggered interaction between bistable units, that is essentially a form of pulse-like coupling, in the presence of noise (stochastic perturbations) and study the escape times dynamics of the nodes in such networks.

\subsection{Sequential escape on networks}

Consider a number of identical and uncoupled systems in quiescent state perturbed by i.i.d. low noise processes. In this case we expect the escapes to occur independently and in random order: the sequence of escapes is a random variable assigning the same probability to each sequence, and the distributions of times of first, second, and $n$th escape will clearly be independent of the sequence. We expect the escapes to occur via visits to partially escaped states of the full system until all have escaped.

In a recent paper \cite{AshCreTsa2017} we highlighted three qualitative regimes that appear for diffusively coupled asymmetric bistable systems as the coupling strength $\beta>0$ is increased: we analyse these effects in detail for a model of epileptic seizure generation in \cite{CreTsaAsh2018}. The regimes can be thought of as emergent phenomena of the noise-perturbed system. For small $\beta$ there is a {\em weak coupling regime} where there is continuation of all partially escaped states and the system is well described by a random sequence that does not necessarily assign the same probability to all sequences - there can be preferred sequences of escape, and the distribution of times of escape may depend on this sequence. For large $\beta$ there is a {\em fast domino regime} of synchronised escapes: the most likely escape path is synchronised and there are no longer any partially escaped attractors. The most interesting, intermediate, range of $\beta$ is the {\em slow domino regime} where some or all of the partially attractors are destroyed but the most likely escape path is not synchronised - this can lead to large but deterministic delays as a domino is deterministically committed to escape but slowly ``topples". As in \cite{AshCreTsa2017} we consider a network of asymmetric bistable systems where the individual systems for $x\in\R$ are of the form
\begin{equation}
\label{eq:onenode}
\dot{x}=f(x,\nu):=-(x-1)(x^2-\nu)
\end{equation}
so that $f=-V'(x)$ and $V(x):=\frac{1}{4}x^4-\frac{1}{3}x^3+\nu(x-\frac{1}{2}x^2)$. For $0<\nu\ll 1$ there is an attractor at $x=x_Q:=-\sqrt{\nu}$ (the ``quiescent'' attractor) and an attractor at $x=x_A:=1$ (the ``active'' attractor) separated by an unstable equilibrium at $x=x_S:=\sqrt{\nu}$, such that the potential has a global minimum $x_A$ and a local minimum at $x_Q$.

In this paper we examine a similar scenario for coupled transient escapes as in \cite{AshCreTsa2017} but with a non-diffusive coupling that is localised in phase space. For this coupling we find no bifurcations of equilibria on increasing $\beta$ but still apparently there is a transition from weak coupling to slow domino regime. Section~\ref{sec:twocoupled} discusses a simple example of two coupled systems of this form and highlights the change in distribution of escape times on increasing $\beta$. In Section~\ref{sec:saddleconnection} we show that this can be explained in terms of a global saddle connection bifurcation.  This shows that destruction of the partially escaped attractors is sufficient, but not necessary, for onset of a slow domino regime. We finish with a discussion in Section~\ref{sec:discuss}.

\section{Local non-diffusive coupling of bistable systems}
\label{sec:twocoupled}

For $i=1,\ldots,N$ we consider a network that evolves according to the It\^{o} stochastic differential equation
\begin{equation}
dx_i=\left[f(x_i,\nu)+\beta\sum_{j\in N_i} h(x_i,x_j)\right]dt + \alpha \,dw_i
\label{eq:Nnode}
\end{equation}
where $N_i$ are the neighbours that provide inputs to node $i$, $\beta$ is the coupling strength, $\alpha$ is the strength of the additive noise and $w_i$ are standard independent Wiener processes. The case studied in \cite{AshCreTsa2017} corresponds to choosing a diffusive coupling function
\begin{equation}
h(x_i,x_j)=(x_j-x_i)
\label{eq:linear}
\end{equation}
where the coupling effect is assumed to be linear with difference in state. In this paper we examine the influence of a non-diffusive coupling of Gaussian form
\begin{equation}
h(x_i,x_j)=H(x_j):=\frac{1}{\sigma\sqrt{\pi}} \exp\left[\frac{(x_j-x_c)^2}{\sigma^2}\right],
\label{eq:gaussian}
\end{equation}
which gives a localised coupling from $x_j$ independent of $x_i$, and acts primarily when $x_j$ is within $\sigma$ of the location of the maximum $x_c$.
In the limit $\sigma\rightarrow 0$, the coupling affects only a small neighbourhood of the lines $\{(x,x_c)\}\cup\{(x_c,x)\}$.

\subsection{Two systems with localised non-diffusive coupling}

\begin{figure}
\begin{center}
\includegraphics[width=12cm]{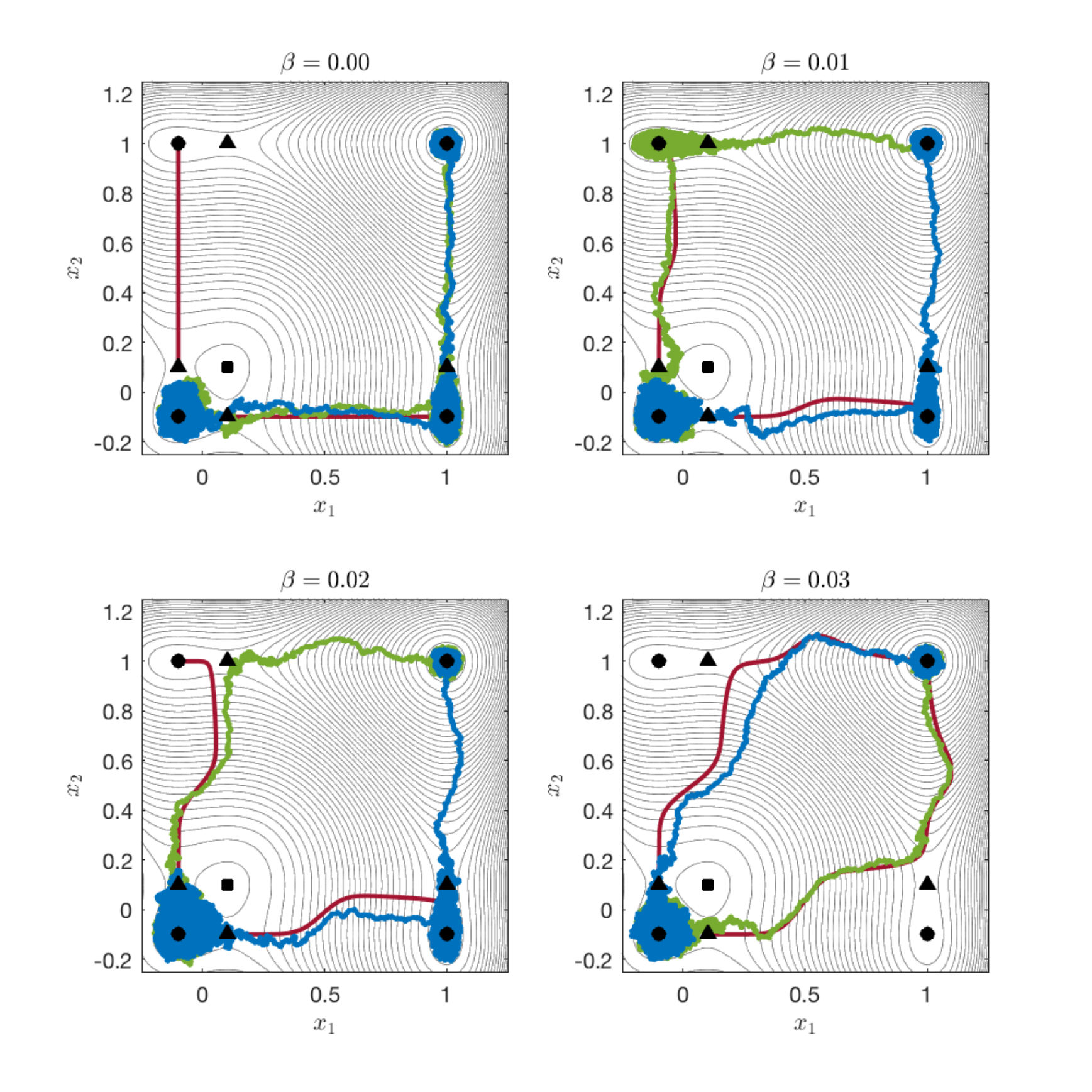}
\end{center}
\caption{Behaviour of the system for varying $\beta$ and $\nu=0.01$. Stable equilibria are denoted $\bullet$, saddles $\blacktriangle$ and repellors $\blacksquare$. The red lines are the unstable manifolds of saddles. Contours for $\beta=0$ show $V(x_1)+V(x_2)$: this is a potential for the uncoupled and noise-free case, and a Lyapunov function for the more general noise- free case $\beta>0$. In each case two realisations are shown (blue and green) for $\alpha=0.02$ that escape from $x_{QQ}$ (bottom left) to $x_{AA}$ (top right). In the bottom left case $\beta=0.02$ close to the saddle connection, note that one trajectory gets trapped in a partially escaped state while the other avoids it.
\label{fig:phaseportraitN2}
}
\end{figure}

Consider the special case of (\ref{eq:Nnode}) for $N=2$ and coupling (\ref{eq:gaussian}), namely
\begin{align}
dx_1&=\left[f(x_1,\nu)+\beta H(x_2)\right]dt + \alpha \,dw_1\nonumber\\
dx_2&=\left[f(x_2,\nu)+\beta H(x_1)\right]dt + \alpha \,dw_2.
\label{eq:2node}
\end{align}
For $\alpha=\beta=0$ the system (\ref{eq:2node}) has stable equilibria at $x_{QQ}=(x_Q,x_Q)$ and $x_{AA}=(x_A,x_A)$ as well as partially escaped stable states at $x_{QA}$ and $x_{AQ}$. Finally there are saddles at $x_{SQ}$, $x_{SA}$, $x_{QS}$, $x_{AS}$ and a source at $x_{SS}$. If we choose default parameters
\begin{equation}
\nu=0.01,~x_c=0.5,~\sigma=0.1
\label{eq:default}
\end{equation}
and vary $\beta$ we find very little change in the bifurcation diagram even for extremely large $\beta>0 $: this is because the coupling terms evaluate to 
$$
\max\{H(x_{Q}),H(x_Q),H(x_S)\}<10^{-6}.
$$
which means $\beta$ must be of order $10^6$ to find a local bifurcation of these states.
However, there is a saddle connection from $x_{QS}$ to $x_{SA}$ for $\beta=\beta_{sc}\approx 0.024$ as can be seen in Fig.~\ref{fig:phaseportraitN2}. This figure illustrates the dynamics for this system, starting at $(x_Q,x_Q)$ for varying values of $\beta$ with noise $\alpha=0.01$. For the noise-free case the unstable manifolds of $x_{QS}$ and $x_{SQ}$ can be seen to lie in the basin of $x_{QA}$ and $x_{AQ}$ for $\beta<\beta_{sc}$ but in the basin of $x_{AA}$ for $\beta>\beta_{sc}$.

\begin{figure}
\begin{center}
\includegraphics[width=12cm]{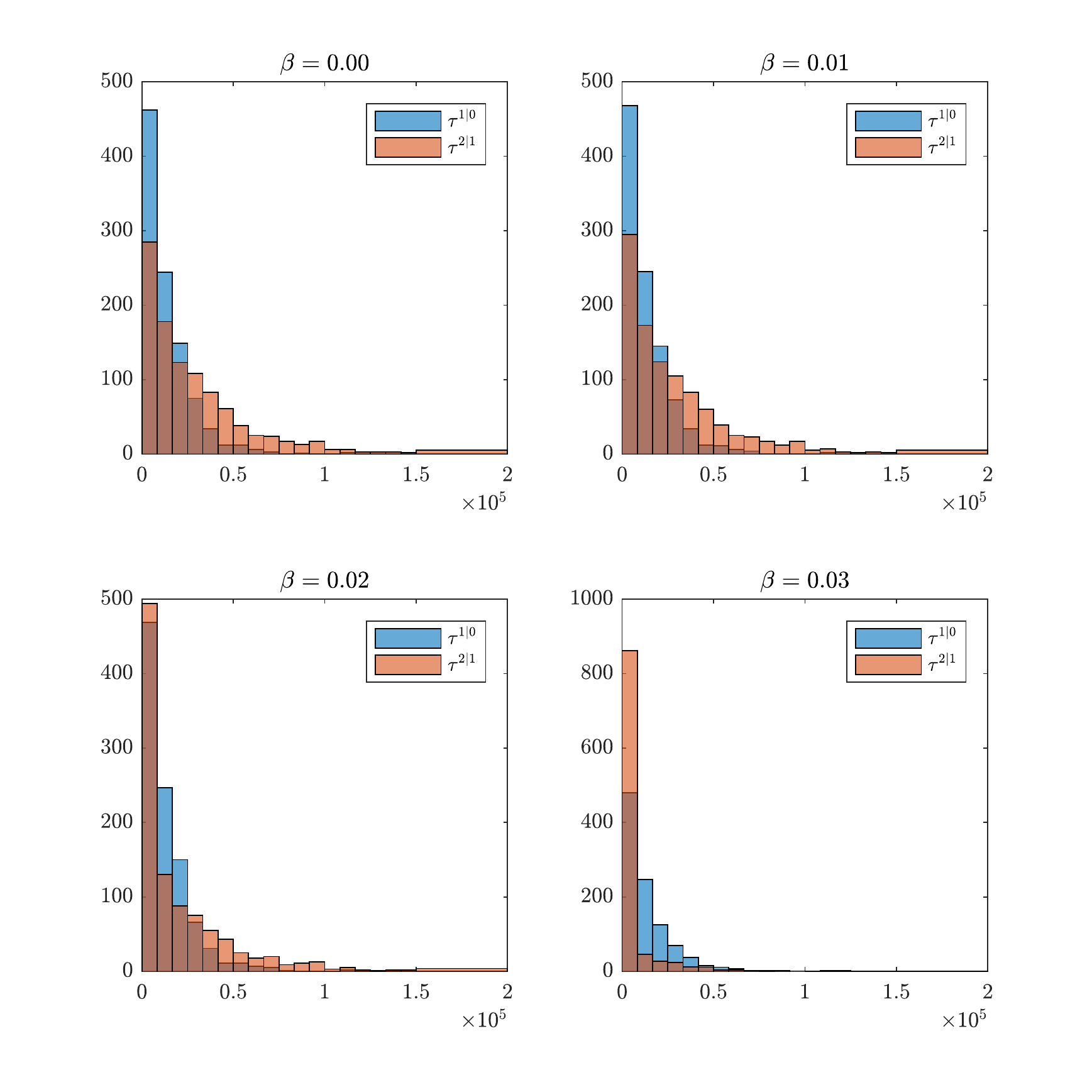}
\end{center}
\caption{Histograms showing changing distributions of first and second escapes on varying $\beta$. The distribution of $\tau^{1|0}$ does not change much with increasing $\beta$; note the scale on the $y$-axis increases with $\beta$. The distribution of $\tau^{2|1}$ changes as $\beta$ increases and for $\beta=0.03$ the second escape time for most simulations is in the smallest bin.
\label{fig:twonodehists}
}
\end{figure}

The sequential escape problem involves global dynamics: as in \cite{AshCreTsa2017,CreTsaAsh2018} we consider the initial state where $x_i=x_{Q}$ for all $i$ and pick a threshold $h$ such that $x_{S}<h<x_{A}$. We define the escape time of the $i$th system to be 
$$
\tau^{(i)}=\inf\{ t>0~:~ x_i(t)\geq h\}.
$$
This is a random variable that depends on system parameters, noise level, network structure and noise realisation. Note that with probability one there is a random permutation (sequence) $s$ such that
$\tau^{(i)}<\tau^{(i+1)}$ for $i=1,\ldots,N-1$. We use this to define the {\em time of $k$th escape} 
$$
\tau^{k}=\tau^{(s(k))}
$$
and $\tau^{k|l}=\tau^k-\tau^l>0$ for any $0\leq l\leq k\leq N$: we define $\tau^0=0$. Fig.~\ref{fig:twonodehists} shows the distributions of the first $\tau^{1|0}$ and second $\tau^{2|1}$ escape times in the two node system.
The histograms are calculated from 1000 realisations of \eqref{eq:2node} with parameters \eqref{eq:default} computed using the stochastic Heun method with timestep $0.001$ with $\alpha=0.02$ and escape threshhold $h=0.8$.
Fig.~\ref{fig:twonodehists} shows a clear change in the distribution of second escape time on varying $\beta$, even though the first escape time is unaffected.

As in \cite{AshCreTsa2017} we define $\E(\tau^{k|l})$ to be the mean of $\tau^{k|l}$, $SD(\tau^{k|l})$ the standard deviation and 
$$
CV(\tau^{k|l})=SD(\tau^{k|l})/\E(\tau^{k|l})
$$
the coefficient of variation. In the case of an exponential (memoryless) distribution of $\tau>0$ note that $CV(\tau)\approx 1$. 

\begin{figure}
\begin{center}
\includegraphics[width=12cm]{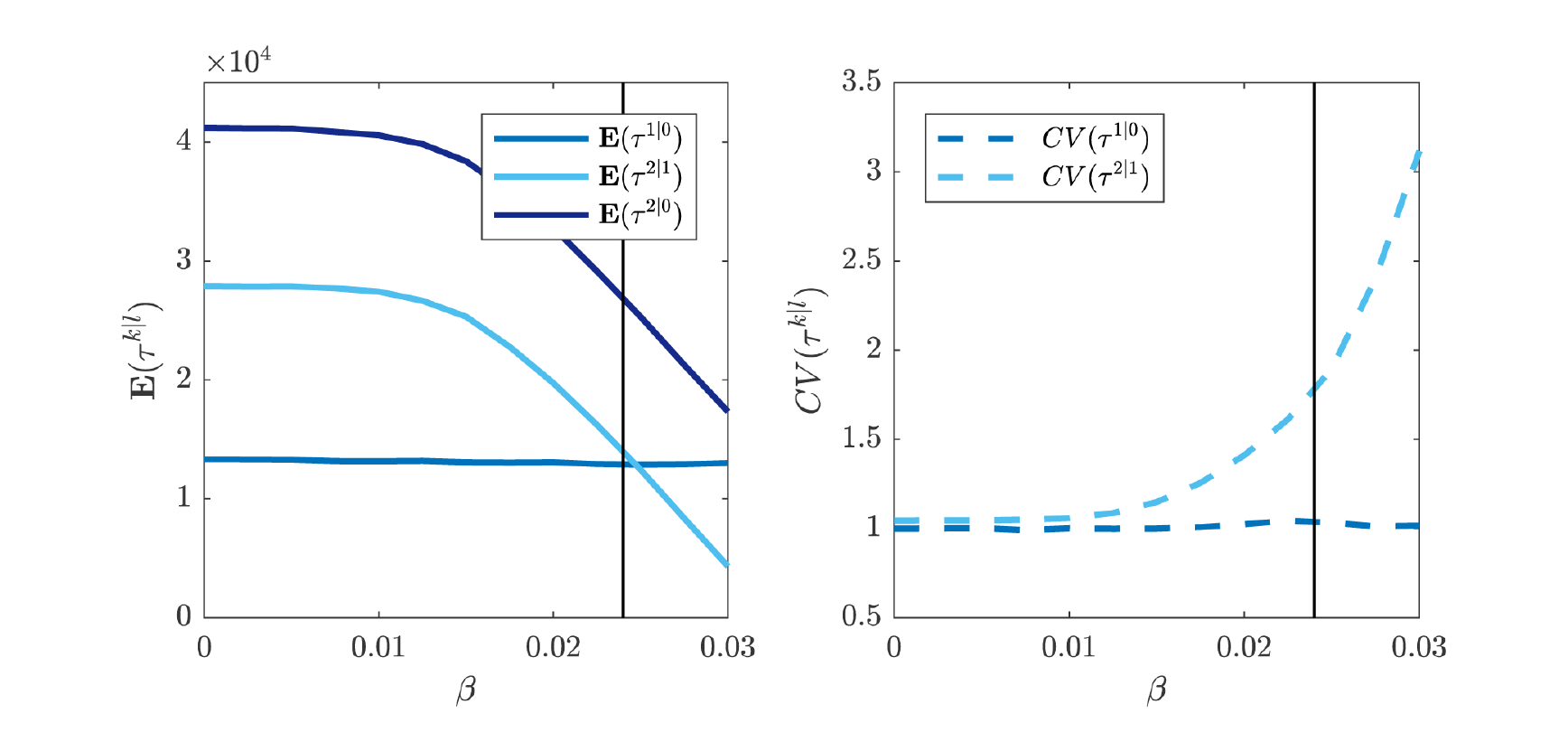}
\end{center}
\caption{Changing mean escape times $\E(\tau^{k|l})$ and $CV(\tau^{k|l})$ on varying $\beta$. For $\beta<0.015$ the second mean escape times $\E(\tau^{2|1}) \approx 2\E(\tau^{1|0})$ , with similar $CV$ values, consistent with uncoupled or weakly coupled behaviour. The mean second escape time $\E(\tau^{2|1}) $ drops quickly and $CV(\tau^{2|1})$ increases exponentially as $\beta$ increases. The mean escape times cross very close to the moment of the saddle connection (black line).
\label{fig:twonodeescapes}
}
\end{figure}

We numerically calculate  $\E(\tau^{k|l})$ and $CV(\tau^{k|l})$  from 1000 realisations computed using the stochastic Heun method as before.
Fig.~\ref{fig:twonodeescapes} shows a clear change in the second escape $\E(\tau^{2|1})$ and $CV(\tau^{2|1})$ on varying $\beta$, and note again that the mean first escape time $\E(\tau^{1|0})$ and $CV(\tau^{1|0})$ are unaffected.

One can understand the emergent properties of the sequential escape probabilities in terms of random Kramers escape from basins of attraction, known to have asymptotic exponential escape rates \cite{FriedlinWentzell1998,BG2006}. This description breaks down for larger coupling in that some escape processes can become predominantly deterministic (corresponding to a small coefficient of variation), or a mixture between processes. 

Friedlin-Wenztell theory \cite{FriedlinWentzell1998} characterises that the mostly likely escape paths from a given stable equilibrium to the next correspond to trajectories that minimise an action integral between one equilibrium and the next: this can be used to define a quasipotential within the basin of the first equilibrium. These trajectories typically first have a diffusion dominated phase where they remain most of the time within the first basin until they find a {\em gate} \cite{Berglund2013}, i.e. a saddle on the basin boundary with the lowest quasipotential. After crossing the gate the trajectory is in a drift dominated phase and will follow a branch of the unstable manifold of the saddle with a rapid motion that gets to the basin of the next equilibrium. 

Note that this description is asymptotic for $\alpha\rightarrow 0$. As discussed in \cite{AshCreTsa2017} the coupling can change the basin of attraction and new gates may bifurcate. In the case of finite $\alpha$ and several possible saddles, the trajectory may escape via a  gate with a higher quasipotential on occasion and this can result in a different sequence of escapes. For finite $\alpha>0$ near the bifurcation this will result in an escape time distribution that is a combination of those for the various possible escape routes.

\section{Saddle connections and onset of weak domino regime}
\label{sec:saddleconnection}

On varying the coupling strength $\beta$ for (\ref{eq:Nnode}) with coupling (\ref{eq:gaussian}) and parameters (\ref{eq:default}) we illustrate the noise-free phase portraits in Fig.~\ref{fig:phaseportraitN2}. Note that, even though there are no apparent changes in the location or stability of the equilibria (or indeed the dynamics associated with the first escape) there is a qualitative change in the unstable manifolds of the equilibria $x_{SQ}$ and $x_{QS}$. Numerically, one can verify that such a saddle connection in the system appears at $\beta=\beta_{sc}\in(0.024,0.025)$, verified by computations similar to those shown in Fig.~\ref{fig:phaseportraitN2}. 

We can find the critical coupling strength $\beta_{sc}$ for (\ref{eq:Nnode}) with coupling (\ref{eq:gaussian}) for the case of the limiting width of interaction $\sigma\rightarrow 0$: note that this corresponds to infinitesimally narrow pulse-like interaction between the units. More precisely, consider a small $d>0$ and suppose that a trajectory passes from $x_2=x_c-d$ to $x_2=x_c+d$ at a point where $\dot{x}_2=\mu>0$ and during this time, approximately $x_2=\mu t$. Hence the main contribution to the change $x_1$ will come from coupling and be approximately
$$
\Delta x_1= \int_{-\infty}^{\infty} \beta\, H(\mu t) \, dt= \frac{\beta}{\mu}.
$$
Hence there will be a saddle connection from $x_{QS}$ to $x_{SA}$ if $\Delta x_1=2\sqrt{\nu}$, i.e. if
$$
\beta_{sc}=2\mu\sqrt{\nu}+O(\sigma)
$$
in the limit $\sigma\rightarrow 0$. Putting in the default parameters (\ref{eq:default}) we have $\mu=f(x_c,\nu)=0.120$ and so asymptotically $\beta_{sc}=0.024$. Comparing to the numerical value approximation, this is in close agreement, even for the moderate value $\sigma=0.1$.

The distribution of second escape times $\tau=\tau^{2|1}$  deviates significantly from memoryless in the case of larger $\beta$: see Fig.~\ref{fig:twonodeescapes}. This is due to there being two very different routes for the second escape. In particular for $\beta$ comparable or larger than $\beta_{sc}$ there will be a probability $0<P<1$ of trajectory after first escape entering the basin of attraction of the partially escaped state $x_{QA}$ or $x_{AQ}$. If it enters this basin, it will escape according to the usual Kramers asymptotics for weak noise. Let $\tau_0$ denote the escape times where it is not captured, and $\tau_1$ the escape times where it is captured. One can compute the mean of the second escape as $\E(\tau)=(1-P)\E(\tau_0)+P\E(\tau_1)$ and its standard deviation $SD(\tau)=\sqrt{(1-P)\E(\tau_0^2)+P\E(\tau_1^2)}$, but note that the coefficient of variation $CV(\tau)= SD(\tau)/\E(\tau)$ has nonlinear dependence on $P$.

\section{Discussion}
\label{sec:discuss}

\begin{figure}
\centering
\includegraphics[width=\linewidth]{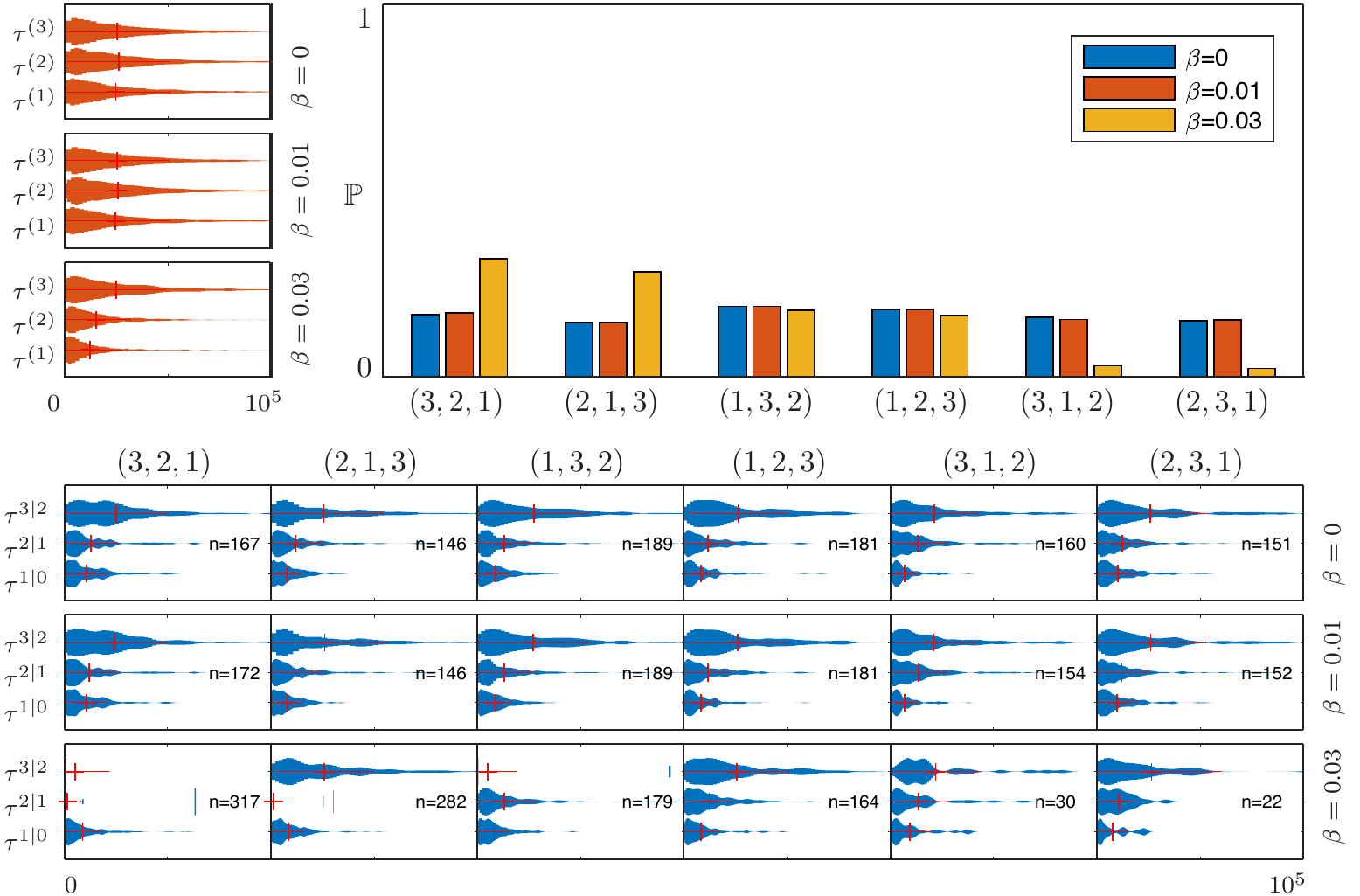}
\caption{Distributions of sequential escapes for the unidirectional chain of three coupled systems (\ref{eq:3node}) with non-diffusive coupling (\ref{eq:gaussian}). Computed using $1000$ realisations; other details as for the two node case. Observe the distribution of the escape times of the nodes (top left) progressively becomes longer for higher $\beta$. The sequential escape time distributions (bottom) are shown as ``violin plots'' where the red vertical bar indicates mean and horizontal bar shows $\pm$ one standard deviation.  The probability of seeing certain sequences (top right) shows that for larger coupling, $\beta=0.03>\beta_c$, the sequences $(3,1,2)$ and $(2,3,1)$ become much less frequent owing to saddle connection bifurcations in the noise-free system.}
\label{fig:threechain}
\end{figure}

The system of two bi-directionally coupled bistable units (\ref{eq:2node}) cannot show emergence of different sequential behaviours of escape, because of permutation symmetry. As a simple example where nontrival sequential behaviour does appear, we briefly consider an analogous system to one considered in \cite{AshCreTsa2017} namely sequential escapes in a uni-directional coupled chain of three units:
\begin{align}
dx_1&=\left[f(x_1,\nu)+\beta H(x_2)\right]dt + \alpha \,dw_1\nonumber\\
dx_2&=\left[f(x_2,\nu)+\beta H(x_3)\right]dt + \alpha \,dw_2\label{eq:3node}\\
dx_3&=\left[f(x_3,\nu)\right]dt + \alpha \,dw_3\nonumber
\end{align}
with $x_i\in\R$, $f$ as in (\ref{eq:onenode}), $H$ as in (\ref{eq:gaussian}) and default parameters (\ref{eq:default}). The coupling parameter $\beta>0$ creates unidirectional forcing on the chain while $\alpha$ modulates the noise level of the standard Wiener processes $w_i$. Note that the system considered in \cite{AshCreTsa2017} had diffusive coupling. In other words, tipping of the third unit affects the second, which in turn affects the first. For the non-diffusive coupling (\ref{eq:gaussian}) we find changes in the relative frequencies of sequences that are quite marked: see Fig.~\ref{fig:threechain}. 

\begin{figure}
\centering
\includegraphics[width=8cm]{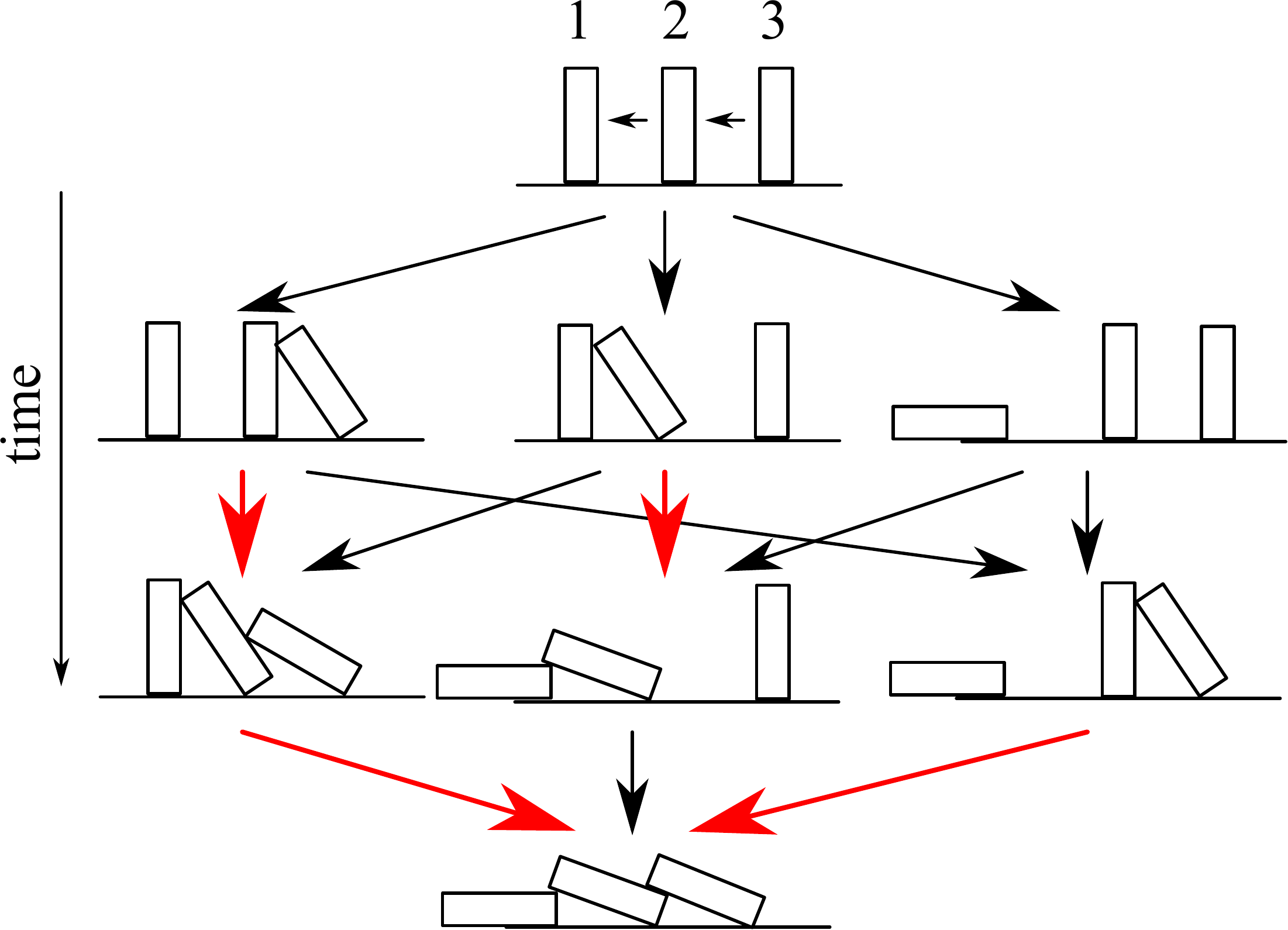}
\caption{Schematic diagram showing possible sequences of escapes for the unidirectional chain of three coupled systems (\ref{eq:3node}) coupled with $3\mapsto 2 \mapsto 1$ and non-diffusive coupling. For the case $\beta=0.03>\beta_c$ there are fast sequence of escapes (shown by red arrows) that greatly reduce the chance of seeing the sequences $(3,1,2)$ and $(2,3,1)$ (see Fig.~\ref{fig:threechain}).}
\label{fig:threechaingraphic}
\end{figure}

The diagram in Fig.~\ref{fig:threechaingraphic} illustrates for $\beta>\beta_c$ the appearance of fast sequences of escapes due to the coupling creating saddle-connections that bypass metastable attractors. Note that since the systems are effectively uncoupled except when one is undergoing a transition, the saddle connections occur at $\beta_c\approx 0.024$ for (\ref{eq:3node}) with standard parameters (\ref{eq:default}). The red arrows on this figure show the fast transitions visible in the sequential escape times shown in Fig.~\ref{fig:threechain}. These rapid escapes greatly reduce the probability of taking one of the less chosen routes $(3,1,2)$ and $(2,3,1)$. In other words, the transition is associated with the most likely escape path from one attractor (corresponding to the unstable manifold of a separating saddle) hits a separating saddle (gate) for a partially escaped attractor.

Varying other coupling parameters we will find other regimes, in particular if $\sigma$ is large or $x_c$ close to one of the equilibrium values $x_{Q,S,A}$ the coupling can lead to bifurcation of equilibria, not just saddle connections. In other words, the transition is associated with the most likely escape path from one attractor (corresponding to the unstable manifold of a separating saddle/gate) hitting the gate for a partially escaped attractor. For example, taking a coupling that mixes both diffusive an localised coupling, we expect that the saddle-node and pitchfork bifurcations noted in \cite{AshCreTsa2017} as organizing the transition from weak coupling to slow and fast domino regimes may interact with the saddle connection in a complex manner. It will be interesting to understand the distributions of escape times that appear in such cases of competition between different escape routes.

As noted previously, sequential escape problems are of relevance to modelling a wide range of problems ranging from epileptogenesis \cite{CreTsaAsh2018} to cell differentiation \cite{Wangetal2010}. In this paper we describe a novel type of emergent behaviour in sequential escapes of coupled systems, associated with a global saddle connection bifurcation. We find changes in the probabilities of seeing certain sequences realised, and changes in the distributions of escape times. It would be good to get a better quantitative understanding of how properties of the escape time distributions (such as coefficient of variation) change on passing through this bifurcation. More specifically, in \cite{AshCreTsa2017} we associated the onset of slow domino regime with the emergence of deterministic escapes and a reduction in CV of the second escape. For the case studied here, the CV apparently increases owing to the mixing of the trapped and non-trapped distributions.

For higher dimensional systems, separating saddles/gates with one-dimensional unstable manifolds should remain of critical importance for understanding the regimes of sequential escape behaviour in the case of varying one coupling parameter. This is because they divide the phase space into basins of attraction from which noise can induce escape. Qualitative regimes of escape from bistable systems with more complex separating invariant sets will also appear at saddle connections between these more complex sets.  

\subsection*{Acknowledgements} 

The authors gratefully acknowledge the financial support of the EPSRC Centre for Predictive Modelling in Healthcare, via grant EP/N014391/1. PA also acknowledges the European Union's Horizon 2020 research and innovation programme for the ITN CRITICS under Grant Agreement number 643073 for providing opportunities to discuss this work with members of the CRITICS network.

\end{document}